\documentclass[letterpaper,12pt]{article}

\usepackage{graphicx}

\usepackage{amssymb}
\usepackage{amsthm}
\usepackage{amsmath}
\usepackage[all]{xy}

\usepackage{enumitem}

\usepackage{centernot}
\usepackage{hyperref}

\usepackage{bm}

\usepackage[cp1251]{inputenc}

\usepackage{tikz-cd}
\usepackage{mathptmx}

\def \Br {{\rm{Br}}}

\def \Pic {{\rm {Pic}}}
\def \val {{\rm {val}}}

\def \Gal {{\rm{Gal}}}

\def \Im {{\rm {Im\,}}}

\def \ev {{\rm{ev}}}

\def \A{{\cal A}}

\def \P{{\mathbb P}}
\def \Spec {{\rm{Spec\,}}}
\def \dim {{\rm{dim\,}}}
\def \ev {{\rm{ev}}}
\def \Hom {{\rm {Hom}}}

\def \Pic {{\rm {Pic}}}

\def\ov{\overline}

\def \Z {{\mathbb Z}}
\def \Q {{\mathbb Q}}

\def\B{{\cal B}}

\def\H{{\rm H}}

\def\O{{\cal O}}
\def\X{{\cal X}}

\def\O{{\cal O}}

\def\sw{{\rm sw}}
\def\fil{{\rm fil}}

\def\beq{\begin{equation} \label}

\newtheorem{theo}{Theorem}
\newtheorem{prop}[theo]{Proposition}
\newtheorem{lem}[theo]{Lemma}
\newtheorem{cor}[theo]{Corollary}
\newtheorem{defi}[theo]{Definition}

\newtheorem{remark}[theo]{Remark}

\newcommand{\bthe}{\begin{theo}}
\newcommand{\ble}{\begin{lem}}
\newcommand{\bpr}{\begin{prop}}
\newcommand{\bco}{\begin{cor}}
\newcommand{\bde}{\begin{defi}}
\newcommand{\ethe}{\end{theo}}
\newcommand{\ele}{\end{lem}}
\newcommand{\epr}{\end{prop}}
\newcommand{\eco}{\end{cor}}
\newcommand{\ede}{\end{defi}}

\date{}
\begin{document}

\title{Evaluation of Brauer elements over local fields}


\author{Evis Ieronymou
}




\maketitle

\begin{abstract}
We study the evaluation maps given by elements of the Brauer group of varieties over local fields.
   We show constancy of the aforementioned maps in several interesting cases.
\end{abstract}

\section{Introduction}
\label{intro}
 Manin, in trying to understand the failure of the Hasse principle combined global class field theory with the Brauer group of a scheme in order to introduce the Brauer-Manin set of a variety over a number field \cite{Man}. This gave birth to the theory of Brauer-Manin obstruction, which has evolved a lot since then and has become an important tool in the study of rational points.

We briefly recall the main points and
 refer the reader to \cite[\textsection 5]{Skbook} for more details on the Brauer-Manin obstruction and its variants or to \cite[\textsection 2]{Wittrc} for a more recent report. Let $X$ be a smooth, proper, geometrically irreducible variety over a number field $k$.
Let $\Br(X)=H^2_{\text{\'{e}t}}(X,\mathbb G_m)$ denote the cohomological Brauer group of $X$. By functoriality of the Brauer group, for any field $F$ containing $k$ we can induce by specialization
 an evaluation map $\ev_{\cal A}:X(F)\to\Br(F)$.
Combining the evaluation maps at all the completions of $k$ 
we get the Brauer-Manin pairing:
$$
\Br(X)\times X(\mathbb A_k)\to \Q/\Z
$$
\noindent where $\mathbb A_k$ is the ring of adeles of $k$. We denote by $X(\mathbb A_k)^{\Br(X)}$ the adelic points othogonal to $\Br(X)$. If we diagonally embed $X(k)$ into $X(\mathbb A_k)$, by the global reciprocity law we have the following chain of inclusions
$$
X(k)\subseteq X(\mathbb A_k)^{\Br(X)} \subseteq X(\mathbb A_k)
$$

 If $X(\mathbb A_k)\neq \emptyset$ and $X(\mathbb A_k)^{\Br(X)}=\emptyset$, we say that there is Brauer-Manin obstruction to the Hasse principle on $X$. If for any $X$ in a class of varieties $\cal C$, we have that $X(\mathbb A_k)^{\Br(X)}\neq\emptyset$ implies $ X(k)\neq\emptyset$, we say that the Brauer-Manin obstruction to the Hasse principle for elements of $\cal C$ is the only one.

Because of conjectures of Bombieri and Lang, it was never expected that the Brauer-Manin obstruction would explain all the failures of the Hasse principle for smooth and proper varieties. Indeed the first example
with $X(\mathbb A_k)^{\Br(X)}\neq\emptyset$ and $ X(k)=\emptyset$ was a bielliptic surface given by Skorobogatov \cite{Skcounter}. Even the more refined \'{e}tale Brauer-Manin obstruction is insufficient to explain all the failures of the Hasse principle as shown by Poonen \cite{Poonen}. The above notwithstanding, we have the following two important conjectures.
The first conjecture was formulated by Colliot-Th\'el\`ene in \cite{CT03}
 and states that if $X$ is rationally
connected, then $X(k)$ is dense in $X(\mathbb A_k)^{\Br(X)}$.
This generalises the same conjecture for geometrically rational surfaces by Colliot-Th\'el\`ene and Sansuc, first asked as a question in \cite{CT.San.con}.
The second conjecture is by Skorobogatov in \cite{Skconj} and states that the Brauer-Manin obstruction to the Hasse principle is the only one for $K3$ surfaces.
Nowadays, there is a large body of literature around the Brauer-Manin obstruction and it has become clear that is important to understand the various evaluation maps given by elements of the Brauer group.

The crucial part of the calculation of the Brauer-Manin obstruction is the calculation of the local evaluation map for each completion of the number field $k$.
In this note we apply deep results of Kato (\cite{Kato1},\cite{Kato2}) to prove the constancy of the local evaluation map in the good reduction case for arbitrary Brauer elements.
This generalises previous known results of Colliot-Th\'{e}l\`{e}ne and Skorobogatov \cite{CTSk} by removing the assumption that the order of the element is coprime to the residual characteristic.
Our approach also recovers many previous results and gives a uniform treatment of all local evaluation maps.
Let $K$ be a finite extension of $\Q_p$ with ring of integers $\O_K$ and let $X$ be a smooth, proper, geometrically irreducible variety over $K$. Assume that there is a model $\X/\O_K$ which is regular with geometrically integral fibres. The novel idea of our approach is the definition of a subgroup $B$ of $\Br(X)$. We show that the relevant evaluation maps are constant for elements of this subgroup.
The definition of $B$ uses Kato's Swan conductor \cite{Kato1}, and depends on the model. 
In some applications we can show that $B$ is the whole Brauer group.
 For example we obtain a relatively simple proof of the following.

{\bf Theorem A}

{\it Let $p$ be an odd prime. Suppose that $\X/\Z_p$ is smooth, proper with geometrically integral fibres, and  either the special fibre is separably rationally connected or the generic fibre is a $K3$ surface. Then

$$
\ev_{\mathcal A}: X(\Q_p)\to \Br(\Q_p)
$$
is a constant map, for any $\mathcal A \in \Br(X)$.}

Note that we actually have more general results in Proposition \ref{regoodRC} for rationally connected varieties and in Proposition \ref{regood} for $K3$ surfaces.
Proposition \ref{regood} already appears in \cite{BN}, cf. Remark 9.

There are some cases where we can replace the assumption on the existence of a good model, with an assumption on the Galois action on $\ell$-adic cohomology.

{\bf Theorem B}

{\it
Let $p$ be an odd prime. Suppose that $X$ is the Kummer surface of an abelian surface over $\Q_p$, and the Galois representation on $\ell$-adic cohomology is unramified for some (any) $\ell\neq p$.
 Then

$$
\ev_{\mathcal A}: X(\Q_p)\to \Br(\Q_p)
$$
is a constant map, for any $\mathcal A \in \Br(X)$.}

We have various results of the above kind in \textsection 4. Note that we formulate the results in the introduction in their simplest form, for ease of exposition. The reader can find more general statements in the main body of the text.

We want to emphasize that besides giving a quick and uniform proof of existing results, our approach has a wide range of applicability.
 For example most of the results of \textsection 5 cannot be obtained by combining existing results in the literature. As a sample, we can show the following (see Remark after Corollary \ref{remmeta})


{\bf Proposition C}

{\it
 Suppose that $X/K$ is a del Pezzo surface which splits over an unramified extension of $K$ and admits a regular proper model with geometrically integral special fibre.
 Then
$$
\ev_{\mathcal A}: X(K)\to \Br(K)
$$
is a constant map, for any $\mathcal A \in \Br(X)$.}

In a different spirit, we also have a result which could be quite useful in computations. Note the interesting feature that the precision needed for the computation depends only on the ground field.

{\bf Proposition D}

{\it
 Suppose that $X/\Q_p$  admits a smooth proper model $\X/\Z_p$ with geometrically integral fibres.
Let $\mathcal A \in \Br(X)[p]$. Then
$$
\ev_{\mathcal A}: X(\Q_p)\to \Br(\Q_p)
$$
factors through
$\X(\mathbb Z_p/p^2)$ for $p$ odd, and through $\X(\mathbb Z_p/p^3)$ for $p$ even.}

Finally, we note that our results have direct implications for a question Swinnerton-Dyer asked the authors of \cite{CTSk}, see \cite[Introduction]{CTSk}. That question was related to his work on density of rational points on certain surfaces, and was the main motivation for \cite{CTSk}.

The outline is as follows. In section 2 we fix notation and state some preliminaries, while in section 3 we prove our main Theorem. Sections 4 and 5 consist of applications in the case of good and bad reduction respectively.

{\it Relation to other work}.
In case that there exists a smooth proper model and $X$ is geometrically simply connected, the constancy of $\ev_{\A}$ for elements of order prime to $p$ follows from \cite[Prop. 2.4]{CTSk}.
Under the same assumptions, for elements of order a power of $p$, one can use \cite[Lem. 7.2, Lem. 7.3]{BN} in order to show constancy of the evaluation map.
Our main result, Theorem \ref{bth} does not use any results from \cite{BN}, but its important consequence, Proposition \ref{simpl.con.l} does (for the case $\ell=p$).

Our approach, on the one hand, gives a short natural proof of our main new result, Theorem A, and, on the other hand, leads to a number of other results including Theorem B and Propositions C and D above (see also Propositions \ref{regoodEnriq}, \ref{hypersurfgood}, \ref{factor.compon} and \ref{curve.res} below).

\section{Preliminaries}
\label{sec:prelim}

We will use the following notation. Given an abelian group $A$ and a positive integer $n$ we denote by $A[n]$ the subgroup of elements annihilated by $n$ and by $A[n^{\infty}]$ the subgroup of elements annihilated by a power of $n$. Given a field $k$, we denote by $\ov k$ a separable algebraic closure of $k$. Given a variety $V$ over $k$, we denote by $\ov V$ or $V_{\ov k}$ the variety $V\otimes_k \ov k$ over $\ov k$.

\begin{itemize}
\item $K$ is a finite extension of $\Q_p$, with ring of integers $\O_K$ and residue field $F$. We denote by $\pi$ a uniformizer of $\O_K$, and by $e$ the absolute ramification index of $K$.

\item
$\mathcal X$ is a faithfully flat, regular, finite type scheme over $\O_K$, with geometrically integral fibres.

\item
$X/K$  is the generic fibre of $\mathcal X/\mathcal{O}_K$.

\item
$Y/F$  is the special fibre of $\mathcal X/\mathcal{O}_K$.

\item $\ell$ is a prime number (the case $\ell=p$ is allowed).

\end{itemize}

We suppose that $X(K)\neq \emptyset$ and we identify $\Br(K)$ with $\Br_0(X)$. We remind the reader that by definition  $\Br_0(X):=\Im(\Br(K)\to \Br(X))$.
\begin{itemize}

\item For a discretely valued field $T$, we denote by  $T_{nr}$ its maximal unramified extension, and we define $\Br_u(T):=\ker(\Br(T)\to \Br(T_{nr}))$.

\item

Let $R=\O_{\X,Y}$, be the local ring of $\X$ at $Y$.  Note that $R$ is a discrete valuation ring. We denote by $R^h$ the henselization of $R$ and by $R^h_{nr}$ the maximal unramified extension of $R^h$.
Hence $R^h_{nr}$ is a strict henselisation of $R$.

\item We set $R^h_{nr,b}:=R^h\otimes_{\O_K}\O_{K_{nr}}$. Equivalently $R^h_{nr,b}$ is the direct limit of the extensions of $R^h$ that correspond to the various extensions $k(Y)/F(Y)$ where $k/F$ is a finite field extension and $k(Y)$ denotes the function field of $Y\otimes_F k$.

\item

We denote by $L$, $L^h$, $L^h_{nr,b}$, $L^h_{nr}$ the fraction field of $R$, $R^h$,  $R^h_{nr,b}$, $R^h_{nr}$ respectively.

\item

We denote by $F(Y)$, $\ov{F}(Y)$ , $\ov{F(Y)}$ the residue field of $R^h$,  $R^h_{nr,b}$, $R^h_{nr}$ respectively.

Note that $F(Y)$ is the function field of $Y$ and $\ov{F}(Y)$ is the function field of $Y\otimes_{F}\ov F$.

\end{itemize}

We will also use some notation from \cite{Kato1}. In particular, see \cite[(1.2)]{Kato1} for the definition of $H^q_n(k)$ and $H^q(k)$, when $k$ is a field. Moreover, see \cite[(1.4)]{Kato1} for the definition of the injective
maps $\lambda_{\pi}:H^{q-1}_n(k)\oplus H^{q}_n(k)\to H^q(K)$ when $K$ is the fraction field of a henselian discrete valuation ring with a chosen uniformiser $\pi$ and residue field $k$ (cf. \cite[Thm 3]{Kato2}).
In this last case Kato also defines an increasing filtration $\{\fil_n H^q(K)\}_{n\geq 0}$ (see \cite[\textsection 2]{Kato1}).
We remind the reader that for any field $k$ we have $H^1(k)=\Hom_{\text{cont}}(\Gal(k^{\text{ab}}/k),\Q/\Z)$ and $H^2(k)=\Br(k)$, see \cite[(1.2)]{Kato1}.

Let $i\geq 1$. We have the following commutative diagram:

\[
\begin{tikzcd}[column sep=tiny]
&&H^1_{\ell^i}(F(Y))\oplus H^2_{\ell^i}(F(Y)) \arrow{r}\arrow{d} &H^1_{\ell^i}(\ov{F}(Y))\oplus H^2_{\ell^i}(\ov{F}(Y))\arrow{r}\arrow{d}&0\arrow{d} \\
\Br(X)\arrow{r}&\Br(L)\arrow{r}&\Br(L^h) \arrow{r}&\Br(L^h_{nr,b})\arrow{r}&\Br(L^h_{nr})
\end{tikzcd}
\]

The vertical maps are injective, and we identify their sources with their images, see \cite[(1.6)]{Kato1} and \cite[Thm 3]{Kato2}.
Note that by \cite[Prop. 6.1 (1)]{Kato1} and \cite[Prop. 1.4.3 (iii)]{CTSkbook} we have that $H^1(F(Y))\oplus \H^2(F(Y))=\Br_u(L^h)$.
Moreover it follows from the definitions in \cite[(1.4)]{Kato1} that $\H^2(F(Y))\subseteq \Br(R^h)$. With a little more work one can show that we actually have equality, but we will not use this fact in the sequel.

For the convenience of the reader we briefly recall some definitions from \cite{Kato1}.
Let $K$ denote a henselian discrete valuation field with valuation ring $\O_K$ and with residue field $F$. Let $\pi$ be a uniformiser of $\O_K$. Kato defines an increasing filtration $\{\fil_n H^q(K)\}_{n\geq 0}$ of the group $H^q(K)$
by
$$
\chi \in \fil_n H^q(K) \Leftrightarrow \{\chi,1+\pi^{n+1}T\}=0 \ \text{in} \ V^{q+1}(\O_K[T])
$$
where $V^{q+1}(\O_K[T])$ denotes a certain direct limit of cohomology groups and $\{, \}$ is induced by cup product.
For the exact definitions we refer the reader to \cite[(1.2), (1.7) and (2.1)]{Kato1}.
The filtration is exhaustive and the swan conductor of $\chi$, denoted by $\sw(\chi)$, is defined as the minimum integer $n$ such that $\chi \in \fil_n H^q(K)$.

\section{The main result}
\label{sec:main.res}
We keep the notation of the previous section.
Before stating the main Theorem we need some preparatory results and definitions. Denote by $\alpha$ the  map
$$
\alpha:\Br(L)\to  \Br(L^h)
$$
and by abuse of notation denote by the same letter the map it induces \\$\Br_u(L)\to  \Br_u(L^h)$.
We have the following lemma

\ble \label{a-1}
$$\alpha^{-1}(H^2_{\ell^i}(F(Y)))\subseteq \Br(R).$$

\ele

\begin{proof}

By \cite[Thm. 3.3]{AB} we have the following commmutative diagram with exact rows

\[
\begin{tikzcd}[column sep=tiny]
0\arrow{r}&\Br(R)\arrow{r}\arrow{d}&\Br_u(L)\arrow{r}\arrow{d}{\alpha}&H^1(F(Y),\Q/\Z)\arrow{r}\arrow{d}{\rm{id}}&0 \\
0\arrow{r}&\Br(R^h)\arrow{r}&\Br_u(L^h)\arrow{r}&H^1(F(Y),\Q/\Z)\arrow{r}&0
\end{tikzcd}
\]

The result follows from the above diagram since $H^2_{\ell^i}(F(Y))\subseteq \Br(R^h)[\ell^i]$.

\end{proof}

We will now define various subgroups of $\Br(X)$.

Let $$ M_{\ell^i}:=\Br(X) \cap \alpha^{-1}(H^1_{\ell^i}(F(Y))\oplus H^2_{\ell^i}(F(Y))).$$

We  have an induced map
$$
M_{\ell^i}\to H^1_{\ell^i}(\ov{F}(Y))\oplus H^2_{\ell^i}(\ov{F}(Y))
$$
and we denote by $\tau_{\ell^i}$ the composition of the map above with the projection to the first factor. Hence we have the following map:
$$
\tau_{\ell^i}:M_{\ell^i}\to H^1_{\ell^i}(\ov{F}(Y)).
$$
In a similar way we define a map:
$$
\tau'_{\ell^i}:M_{\ell^i}\to H^1_{\ell^i}(F(Y)).
$$

We now define ${}_{\ell}B\subseteq {}_{\ell}M\subseteq \Br(X)$ as follows:

$${}_{\ell}M:=\bigcup_{i\geq 0} M_{\ell^i}.$$
$${}_{\ell}B:=\bigcup_{i\geq 0} \ker(\tau_{\ell^i}).$$

We caution the reader that what Kato denotes by $K$ and $F$ in our case are $L^h$ and $F(Y)$ respectively.
The filtration of Kato on $\Br(L^h)$  (see \cite[\textsection 2]{Kato1}) induces a filtration on $\Br(X)$ and $\Br(L)$.
With respect to this filtration we have the following important result.
\bpr \label{dias}
\begin{enumerate}
\item
 ${}_{\ell}M[\ell^{\infty}]=\fil_0(\Br(X))[\ell^{\infty}]$.
\item If $\ell\neq p$, then ${}_{\ell}M[\ell^{\infty}]=\Br(X)[\ell^{\infty}]$.
\end{enumerate}
\epr
\begin{proof}
This follows from \cite[Prop. 6.1 (1)]{Kato1} cf. last paragraphs of Preliminaries.

\end{proof}
Denote by $B$ the subgroup of $\Br(X)$ generated by the ${}_{\ell}B$ as $\ell$ varies through all the primes. Note that $B$ depends on the regular model with geometrically integral fibres $\X/\O_K$.
When we want to be explicit about this depedence we will denote $B$ by ${}^{\X}\Br(X)$.
We can now state and prove our main result.

\bthe
  \label{bth}
Let $\mathcal A \in B$. Then

$$
\ev_{\mathcal A}: \X(\O_K)\to \Br(K)
$$
is a constant map.

\ethe

\begin{proof}
 We may assume that $\A \in {}_{\ell}M$, for some $\ell$. Fix $i$ with $\mathcal A \in \ker(\tau_{\ell^i})$. Let $\chi\in H^1_{\ell^i}(F)$ denote a character of the Galois group of $F$ of order $\ell^i$ and denote by $\phi:\H^1_{\ell^i}(F) \to \H^1_{\ell^i}(F(Y))$ the natural map. It follows from the definitions that the image of $\chi$ in $\Br(K)$ is the cyclic algebra $\B=(\chi,\pi)$, and $\alpha(\B)$ is the image of $\phi(\chi)$ in $\Br(L^h)$.
Since $Y$ is geometrically integral, $\phi$ is injective and hence $\alpha(\B)$ has order $\ell^i$.
Now, by the inflation-restriction sequence we have
$$
\ker(H^1_{\ell^i}(F(Y))\to H^1_{\ell^i}(\ov F(Y))=H^1(G,\Z/\ell^i)
$$
where $G=\Gal(\ov{F}(Y)/F(Y))\cong \hat \Z$. Therefore
$\ker(H^1_{\ell^i}(F(Y))\to H^1_{\ell^i}(\ov F(Y))$ is the subgroup generated by $\alpha(\B)$. As $\tau'_{\ell^i}(\A)\in \ker(H^1_{\ell^i}(F(Y))\to H^1_{\ell^i}(\ov F(Y))$ by assumption, we deduce that $\A-n\B\in \alpha^{-1}(H^2_{\ell^i}(F(Y)))$
for some $n\in \Z$. By Lemma \ref{a-1}, we have that $\A-n\B\in \Br(R)$. Hence  $\A-n\B\in \Br(R)\cap \Br(X)=\Br(\X)$  where the last equality follows from \cite[Thm. 3.7.6]{CTSkbook}.
Since $\ev_{{\cal C}}$ is identically zero for ${\cal C} \in \Br(\X)$ it follows that $\ev_{\mathcal A}$ sends everything to $n\B$.

\end{proof}

To state our next result and make the link to \cite{CTSk} more explicit, denote by $\beta$ the  map
$$
\beta:\Br(L)\to  \Br(L^h_{nr,b}).
$$

\bco
 \label{non.p.result}

Let $\A \in \Br(X)$. Suppose that $\A\in \ker(\beta)$. Then
$$
\ev_{\mathcal A}: \X(\O_K)\to \Br(K)
$$
is constant

\eco

\begin{proof}
We can suppose that the order of $\A$ is a power of $\ell$. By \cite[Prop. 6.1 (1)]{Kato1}, we have that $\A\in \fil_0(\Br(X))$. Therefore we have that $\A \in {}_{\ell}M$ by Proposition $\ref{dias}$.
 It is then clear from the definitions that $\A\in {}_{\ell}B$, and so we conclude by Theorem \ref{bth}.

\end{proof}

\begin{remark}
  Note that $\ker(\Br(X)\to \Br(X\otimes_K K_{nr}))\subseteq \ker(\beta)$ and so this corollary recovers \cite[Lemma 2.2 (ii)]{CTSk}, where they additionally assume that $\X/\O_K$ is proper.
Note also that in the proof of \cite[Prop. 2.3]{CTSk} it is shown that if $H^1(X,\O_X)=0$ and the Neron-Severi group $NS(X\otimes_K \ov K)$ is torsion-free then $\ker(\Br(X)\to \Br(X\otimes_K K_{nr}))=\Br_1(X)$.
\end{remark}

We can combine Theorem \ref{bth} with two results from \cite{BN} to obtain the following.
\bpr
  \label{simpl.con.l}
Suppose that $\X/\O_K$ is smooth and that the maximal pro-$\ell$ quotient of the geometric fundamental group of $Y$ is trivial.
Let $\mathcal A \in \Br(X)$ and suppose one of the following:

\begin{enumerate}[label=(\roman*)]
\item
$\mathcal A \in {}_{\ell}M$;
\item $\mathcal A \in \Br(X)[\ell^{\infty}]$ and $\fil_0(\Br(X))[\ell^{\infty}]=\Br(X)[\ell^{\infty}]$;
\item $\mathcal A \in \Br(X)[\ell^{\infty}]$ and $\ell \neq p$;
\item $\mathcal A \in \Br(X)[\ell^{\infty}]$, $\ell = p$, $H^0(Y,\Omega_Y^1)=0$ and $e<p-1$.
\end{enumerate}

 Then

$$
\ev_{\mathcal A}: \X(\O_K)\to \Br(K)
$$
is a constant map.

\epr

\begin{proof} Cases (ii) and (iii) follow from Proposition \ref{dias} and case (i).
Case (iv) follows from \cite[Lemma 7.2]{BN}
 and case (ii). It remains to prove case (i).
Fix $i$ with $\mathcal A \in M_{\ell^i}$. By Proposition  \cite[Prop. 3.1]{BN} for the case $\ell=p$ and by   \cite[proof of Prop. 2.4]{CTSk} for the case $\ell\neq p$, we see that $\tau_{\ell^i}(\A)$ lies in $H^1(Y\otimes_F \ov F, \Z/\ell^i)$, which is trivial by assumption.
Therefore $\A\in {}_{\ell}B$ and we are done by Theorem \ref{bth}.

\end{proof}

\section{Applications to good reduction}
\label{sec:good.red}

We keep the notation and assumptions of the previous section.
In this section we will use some more notions from  \cite{Kato1}. In particular see \cite[pg. 121]{Kato1} for the definition of $\sw_{\mathfrak p}(\A)$,
when $S$ is a normal irreducible scheme with function field $L$, $\mathfrak p\in S^1$ (i.e. $\mathfrak p\in S$ and $\dim(\O_{S,\mathfrak p})=1$), and $\A\in \Br(L)$.

Our first result concerns rationally connected varieties. In respect to the relation of the hypothesis to the Brauer-Manin obstruction we note the following which follows from \cite[IV Thm 3.11]{Kollar}: if $k$ is a number field and $X/k$ is smooth, proper and rationally connected, then there is a finite set of places $T$ containing the archimedean places such that $X\otimes_k{k_v}$ admits a smooth, proper model with separably rationally connected special fibre for any $v\notin T$.

\bpr

\label{regoodRC}
Suppose that $\X/\O_K$ is smooth and proper, and the special fibre is separably rationally connected.
Then

$$
\ev_{\mathcal A}: X(K)\to \Br(K)
$$
is a constant map, for any $\mathcal A \in \Br(X)$.

\epr

\begin{proof}
We can assume that the order of $\A$ is a power of a prime $\ell$. The special fibre is geometrically simply connected, see e.g. \cite[Cor. 3.6]{Deb}.
Since $\X/\O_K$ is proper, we can conclude by Proposition \ref{simpl.con.l}, once we show that $\fil_0(\Br(X))[\ell^{\infty}]=\Br(X)[\ell^{\infty}]$ in the case $\ell=p$.
 Note that by \cite[Cor. IV. 3.8]{Kollar} and its proof we have that $H^0(Y,\Omega_Y^1)=H^0(Y,\Omega_Y^2)=0$ . Hence it follows from  \cite[Lemma 7.1]{BN} that $\A\in \fil_0(\Br(X))$. 

\end{proof}

Our next result concerns $K3$ surfaces.
\bpr
\label{regood}
Suppose that $\X/\O_K$ is smooth and proper, $X/K$ is a $K3$ surface and $e<p-1$.
Then

$$
\ev_{\mathcal A}: X(K)\to \Br(K)
$$
is a constant map, for any $\mathcal A \in \Br(X)$.
\epr

\begin{proof}
We can assume that the order of $\A$ is a power of a prime $\ell$. Under our assumptions $\Pic(\X)=\Pic(X)$ and $\Omega_{\X/\O_K}$ is locally free. These together with the fact   that $\omega_X$ is trivial and that the formation of exterior powers of the sheaf of relative differentials commutes with base change, imply that $\Omega^2_{\X/\O_K}$ is also trivial. Hence $\omega_Y$ is trivial as well. Moreover, since $\X/\O_K$ is smooth and proper it follows from
\cite[Cor. VI. 4.2]{Milne} that the Betti numbers of the generic and the special fibre are equal. Hence $Y$ is a $K3$ surface by the classification of surfaces in characteristic $p$, see e.g. \cite{Bomb}. Therefore $H^0(Y,\Omega_Y^1)=0$ by a Theorem of Rudakov
and \u Safarevi\u c \cite{ShRu}, see also \cite{Lang.Nyg}. Hence the result follows from Proposition \ref{simpl.con.l}, since $\X/\O_K$ is proper.

\end{proof}

\begin{remark} The above result is already in \cite[Remark 7.5]{BN}. They showed it by combining some of their results with results from \cite{CTSk} as explained in the {\it Relation to other work} paragraph of the introduction.
\end{remark}

The next corollary is interesting as it does not explicitly mention the existence of a good model in its assumptions.

\bco
 \label{regoodKumsurf}

Suppose that $X$ is the Kummer surface of an abelian surface $A$ over $K$, the $\Gal(\ov K/K)$-representation $H^2(X_{\ov K},\Q_{\ell})$ is unramified for some (any) $\ell\neq p$ and $e<p-1$. Then

$$
\ev_{\mathcal A}: X(K)\to \Br(K)
$$
is a constant map, for any $\mathcal A \in \Br(X)$.

\eco

\begin{proof}
By Proposition \ref{regood}, it suffices to show the existence of $\X/\O_K$ which is smooth, proper and has $X$ as generic fibre. This follows from \cite[Thm. 4.1]{Mats}.

\end{proof}

We record a result for Enriques surfaces.
\bpr
\label{regoodEnriq}

Suppose that $\X/\O_K$ is smooth and proper and that the generic fibre is an Enriques surface.
Let $\mathcal A \in \Br(X)$ have odd order. Then

$$
\ev_{\mathcal A}: X(K)\to \Br(K)
$$
is a constant map.

\epr

\begin{proof}
It is well known that the geometric fundamental group of the generic fibre is isomorphic to $\mathbb Z/2$, see e.g. \cite[Lemma VIII.15.1]{CCSurf} .
By \cite[Thm. X. 3.8]{SGA 1} the geometric fundamental group of the special fibre is a quotient of $\mathbb Z/2$.
 Under our assumptions $\Pic(\X)=\Pic(X)$ and $\omega_X^2$ is trivial. By the same reasoning as in the proof of Proposition \ref{regood} it follows that $\omega_Y^2$ is trivial. Moreover the Betti numbers of the generic and the special fibre are equal, since $\X/\O_K$ is smooth and proper \cite[Cor. VI. 4.2]{Milne}. Hence $Y$ is an Enriques surface surface by the classification of surfaces in characteristic $p$, see e.g. \cite{Bomb}.
 We can assume that the order of $\A$ is a power of an odd prime $\ell$.
Since $\X/\O_K$ is proper, we can conclude by Proposition \ref{simpl.con.l}, once we show that $\fil_0(\Br(X))[\ell^{\infty}]=\Br(X)[\ell^{\infty}]$ in
the case $\ell=p$. We have that $H^0(Y,\Omega_Y^1)=H^0(Y,\Omega_Y^2)=0$ by \cite[Prop. I. 1.4.1]{Dolg}.
Hence it follows from \cite[Lemma 7.1]{BN} that $\A\in \fil_0(\Br(X))$.

\end{proof}

\begin{remark} If $X$ is an Enriques surface then $\Br(\ov X)\cong \Z/2$, and so an element of odd order in $\Br(X)$ is necessarily algebraic. I do not know of examples where $H^1(K,\Pic \ov X)\cong \Br_1(X)/\Br_0(X)$ has odd torsion.

\end{remark}


We have the following for smooth surfaces in $\mathbb P^3$.

\bpr
\label{hypersurfgood}

Suppose that $\X/\O_K$ is smooth and proper. Assume that the special fibre is isomorphic to a surface in $\mathbb P^3$ and $e<p-1$.
 Then

$$
ev_{\mathcal A}: X(K)\to \Br(K)
$$
is a constant map, for any $\mathcal A \in \Br(X)$.

\epr

\begin{proof}
The special fibre is geometrically simply connected and it follows from Bott vanishing for $\P^3$  that $H^0(Y,\Omega_Y)=0$. Therefore we are done by Proposition \ref{simpl.con.l}.

\end{proof}

The assumption for the crystalline cohomology in the next result, is not easy to check in general.
However we note that we do get a new proof of the constancy of the evaluation maps for some classes of $K3$ surfaces, since such surfaces are known to have torsion free crystalline cohomology.

\bpr\label{regoodKumsurf}

Assume that $K=\Q_p$ and $p\neq 2$.

Suppose that $\X/\O_K$ is smooth and proper, and the generic fibre $X/K$ is one of the following:

\begin{enumerate}
\item
 geometrically Kummer;
\item
 smooth complete intersection of dimension at least 2;
\item
isomorphic to an $n$-dimensional subvariety of $\mathbb P^r$, with $2n-r\geq 1$.

\end{enumerate}

Assume moreover that  $H^1_{\rm{cris}}(Y/W)$ and $H^2_{\rm{cris}}(Y/W)$ are torsion-free.

Let $\mathcal A \in \Br(X)$. Then

$$
ev_{\mathcal A}: X(K)\to \Br(K)
$$
is a constant map.

\epr

\begin{proof} Take an embeding $K\to \mathbb C$ and let $X_{\mathbb C}=X\otimes_K \mathbb C$. Then it is well known that $X_{\mathbb C}$ is simply connected and has first Betti number equal to 0, see e.g. \cite[Thm. 1 and Thm. 2]{Span} and \cite[Thm. 3.2.1 and Cor. 3.2.2]{Laz}. This implies that $X$ is geometrically simply connected (\cite[Cor. X. 1.8]{SGA 1}) and that $H^1_{dR}(X)=0$.
The fact that $H^1_{dR}(X)=0$ and our assumptions on the crystalline cohomology imply that $H^1_{dR}(Y)=0$ (see \cite[\textsection 1.3]{Ill1}).
Therefore it follows from the Hodge spectral sequence for $Y$ and  \cite[Cor. 2.5]{Del} that $H^0(Y,\Omega_Y)=0$. We can conclude by Proposition \ref{simpl.con.l}.

\end{proof}

\begin{remark} If $X$ is a smooth complete intersection of dimension at least 3, then $\Br(X)=\Br_0(X)$ by \cite[Prop. A.1]{CT}.
\end{remark}

 Before we state the next Proposition, we need some preparation.
 Suppose that $\X/\O_K$ is smooth and let $\mathfrak p\in \X$ be the generic point of $Y$. Let $P_0\in \X(F)\subset \X$
 and let $A=\O_{\X,P_0}$.
 We will need the following lemma (cf. \cite[\textsection 3.2]{BLR}).

\ble\label{blup}

Let $\tilde{\X}\to \X$ denote the blow-up of $\X$ at $P_0$.
Let $\cal Z$ be $\tilde {\X}$ minus the strict transform of $Y$.
Then

\begin{itemize}
  \item[(i)] ${\cal Z}/\O  _ K$ is smooth;
  \item[(ii)] the map ${\cal Z}(\O_K)\to \X(\O_K)$ surjects to the elements of $\X(\O_K)$ that reduce to $P_0$;
  \item[(iii)] 
   for any integer $i\geq 1$, the map ${\cal Z}(\O_K)\to {\cal Z}(\O_K/\pi^{i})$  factors through the map ${\cal Z}(\O_K)\to \X(\O_K/\pi^{i+1})$.

\end{itemize}

\ele

\begin{proof}
 It follows from \cite[Thm. 3.1]{Valla}, \cite[Prop. 6.1.5]{EGAIV2} and \cite[Lemma 01V8]{Stacks} that ${\cal Z}/\O  _ K$ is smooth.

  Let $s\in \X(\O_K)$ be an element that reduces to $P_0$. It is easy to see that $s$ lifts uniquely to an element of $\tilde{\X}(\O_K)$ by the universal propery of blowing-up, see e.g. \cite[Prop. II. 7.14]{Har}.
 In order to show (ii) and (iii) we can work locally and replace $\X$ by $\Spec(A)$.
 We will follow the proof of \cite[Prop. II. 7.14]{Har}.
 Let $\pi,a_1,\cdots,a_m$ be a regular system of parameters of $A$. These define a closed embedding of $\tilde{\X}$ in $\mathbb P_{A}^m$.
 The exceptional divisor $E$ is isomorphic to $\mathbb P_{F}^m$ and the intersection of $E$ with the strict transform of $Y$ corresponds to the hyperplane in $E$ with first coordinate equal to $0$.
   Following through the definitions we see that $s$ lifts to an element of ${\cal Z}(\O_K)$ which has coordinates $\left( \frac{s(a_1)}{\pi},\cdots,\frac{s(a_n)}{\pi} \right )$ in $\mathbb A_{A}^m$.
  We can deduce parts (ii) and (iii) from this explicit description.

   Let's show part (iii). All morphisms that follow are ring homomorphisms.  Note that $s$ is a section of the structure map $\O_K \to A $, and that $\pi,a_1,\cdots,a_m$ generate $\frak m$, the maximal ideal of $A$.

  We have
  $$
  \O_K\to A\xrightarrow{s} \O_K,
  $$
  that  $F=\O_K/\pi=A/\frak m$, $s(\pi)=\pi$ and $\val( s(a_i))\geq 1$.

  Let $\tilde s,\tilde \tau \in \mathcal Z(\O_K)$ denote the lifts of $s,\tau \in \X(\O_K)$, where $\tau$ also reduces to $P_0$.
  We have to show the following: If $s,\tau$ have the same image in $\X(\O_K/\pi^{i+1})$ then $\tilde s,\tilde \tau$ have the same image in $\mathcal Z(\O_K/\pi^{i})$.
  From the explicit description of the lifts, this readily translates to showing the following
  $$
  s(a_j)\equiv \tau(a_j) \mod \pi^{i+1} \Rightarrow \frac{s(a_j)}{\pi}\equiv \frac{\tau(a_j)}{\pi} \mod \pi^{i}
  $$
  for all $1\leq j\leq m$. It is easy to see that the last displayed implication holds.



\end{proof}

We now recollect some things from \cite[\textsection 7 and \textsection 8]{Kato1}. Let $A=\O_{\X,P_0}$ as before.
We use the notation of \cite[\textsection 7]{Kato1}.
 In that notation $k=F$ and we caution the reader that the $K$ appearing there denotes the function field of $\X$ and so it is different from the $K$ in this paper. We will need the special case $q=2$, $r=1$, $\mathfrak p=\mathfrak p_1=\pi A$.
 Note that an element $\mathcal A\in \Br(X)$ satisfies the assumptions at the beginning of \cite[\textsection 7]{Kato1}.
Moreover since $F$ is a finite field we have that $\Omega^i_F$ is trivial for $i>0$. Therefore if $\sw_{\mathfrak p}(\A)\geq 1$ then $\A$ is not strongly clean with respect to $A$
(see \cite[Def. 7.4]{Kato1} for the definition of strongly clean).


We can now state and prove our next result, which can be useful in actual computations. Note that the $n$ appearing in the statement of the Proposition depends only on $K$.
\bpr\label{precision}
Suppose that $\X/\O_K$ is smooth. Let $n$ be the smallest integer that is greater than $\frac{e}{p-1}$.
Let $\mathcal A \in \Br(X)[p^t]$. Then

$$
\ev_{\mathcal A}: \X(\O_K)\to \Br(K)
$$
  factors through $\X(\O_K/\pi^{n+te})$.
\epr
\begin{proof}
Let $\mathfrak p\in \X$ be the generic point of $Y$. If $\sw_{\mathfrak p}(\A)=0$ then by \cite[Prop. 3.1]{BN} $\ev_{\A}$ factors through $\X(\O_K/\pi)$ and we are done. 
Hence we suppose that $\sw_{\mathfrak p}(\A)\geq 1$. Let $P_0\in \X(F)\subset \X$ and let $A=\O_{\X,P_0}$. Let $\tilde{\X}\to \X$ denote the blow-up of $\X$ at $P_0$.
The residue field of $A$ is a finite field and as explained before the Proposition it follows that $\A$ is not strongly clean with respect to $A$.
 Therefore if $\nu\in \tilde \X$ denotes the generic point of the exceptional divisor, it follows from \cite[Thm 8.1]{Kato1} that $\sw_{\nu}(\A)<\sw_{\mathfrak p}(\A)$. Let $\cal Z$ be $\tilde {\X}$ minus the strict transform of $Y$.

 If $\sw_{\nu}(\A)=0$ then by \cite[Prop. 3.1]{BN} applied to $\cal Z$, the map $ev_{\A}$ is constant on the elements of ${\cal Z}(\O_K)$
with the same image in ${\cal Z}(\O_K/\pi)$.
Moreover, it follows from Lemma \ref{blup} that if $s_1,s_2\in \mathcal X(\O_K)$ reduce to $P_0$ and have the same image in $\X(\O_K/\pi^2)$, then their lifts
 in $\mathcal Z(\O_K)$ have the same image in $\mathcal Z(\O_K/\pi)$. Therefore $\ev_{\A}$ is constant on the elements of $\X(\O_K)$ that reduce to $P_0$ and have the same image in $\X(\O_K/\pi^2)$.
If $\sw_{\nu}(\A)\neq 0$ then we repeat the argument starting with a point in ${\cal Z}(F)$. If this process stops after $n+te-1$ times then the Proposition follows in a similar vein to the above, by using Lemma \ref{blup} repeatedly (in the end we apply the same argument for each point in $\X(F)$) .
Therefore it suffices to show that $\sw_{\mathfrak p}(\A)\leq n+te-1$.

By \cite[proof of Lemma 2.4]{Kato1}, in particular by \cite[(2.4.1)]{Kato1}, it follows by an easy induction that
every element of $1+\pi^{n+ke}R$ is a $p^k$-power in $R^*$ for every integer $k\geq 1$, where $R$ is the ring denoted $(\O_K[T])^{(h)}$ in \cite{Kato1}.
We remind the reader that Kato's $K$ in \cite{Kato1} is $L^h$ in our notation. In particular $1+\pi^{n+te}T$ is a $p^t$-power.
Since the cup product is bilinear and $p^t\mathcal A=0$ it follows from \cite[Def. 2.1 and Def. 2.3]{Kato1} that $\sw_{\mathfrak p}(\A)\leq n+te-1$.

\end{proof}

\begin{remark} In contrast, if the order of $\A$ is coprime to $p$ then it is well-known that $\ev_{\A}$ factors through $\X(\O_K/\pi)$, see e.g. \cite[proof of Prop. 2.4]{CTSk}.
\end{remark}

\section{Applications to bad reduction}
\label{sec:bad.red}

In this section we assume that $\X/\O_K$ is regular and $X/K$ is smooth, projective and geometrically integral. Let $Y=\X\otimes F$ be the special fibre. {\it We assume that the  irreducible components of $Y$ are geometrically integral.} Note that an element of $\X(\O_K)$ will intersect the smooth locus of a unique irreducible component of $Y$.
We thus have a map $\X(\O_K)\to \rm{IrredComp}(Y)$.
Our first result in this section is also an illustration of the robustness of our approach.


\bpr\label{factor.compon}

Assume  that $\A \in \ker(\Br(X)\to \Br(X\otimes_K K_{nr}))$. Then


$$
\ev_{\mathcal A}: \X(\O_K)\to \Br(K)
$$
factors through $\rm{IrredComp}(Y)$
\epr

\begin{proof}

Let $$\rm{IrredComp}(Y)=\{Y_1,\cdots,Y_n\}.$$ Fix $1\leq i \leq n$ and set $\X'=\X-\bigcup_{j\neq i} Y_j$. By Corollary \ref{non.p.result} and the remark below it, we see that the map
$\X'(\O_K)\to \Br(K)$ given by evaluating at $\A$ is constant. This is what we had to prove. 

\end{proof}

\begin{remark} The above gives a new proof of part of the main result in \cite{Bri.eff}, see \cite[Thm. 1]{Bri.eff}
\end{remark}

\bco
\label{remmeta}

Let $I=\Gal(\ov K/ K_{nr})$ and assume  that $H^1(I,\Pic(\ov X))=0$.\\ Let $\mathcal A \in \Br_1(X)$. Then


$$
\ev_{\mathcal A}: \X(\O_K)\to \Br(K)
$$
factors through $\rm{IrredComp}(Y)$.
\eco
\begin{proof} The condition $H^1(I,\Pic(\ov X))=0$ and the fact that $I$ has cohomological dimension $1$, imply via the Hochschild-Serre spectral sequence that $\Br_1(X\otimes_K K_{nr})$ is trivial.
 Therefore it follows immediately that $$\Br_1(X)=\ker(\Br(X)\to \Br(X\otimes_K K_{nr}))$$
and we can conclude by Proposition \ref{factor.compon}.

\end{proof}

\begin{remark} For example let $X$ be a del Pezzo surface which admits a regular proper model such that the special fibre is geometrically integral. If $X$ splits over an unramified extension of $K$ or more generally if $\Br(X\otimes_K K_{nr})$ is trivial then the  map $$\ev_{\A}:X(K)\to \Br(K)$$ is constant for any $\A\in \Br(X)$. See \cite[\textsection 4]{Br2} for conditions on the special fibre that ensure that $\Br(X\otimes_K K_{nr})$ is trivial.\end{remark}

We also record the following result for curves, which is potentially useful for computations in specific cases. We note that the assumptions are the necessary ones for the argument to go through.
\bpr\label{curve.res}
Suppose that $\X/\O_K$ is regular and proper, $X/K$ is a smooth, geometrically integral curve of genus $g$ and $Y\otimes_F \ov F$ is irreducible with a unique singular point $P$.
We assume that there is only one tangent direction to $\ov Y$ at $P$, and that $P$ has multiplicity $r$ with $g=\frac{r(r-1)}{2}$.

 Let $\mathcal A \in \Br(X)$ have order coprime to $p$. Then

$$
\ev_{\mathcal A}: X(K)\to \Br(K)
$$
is a constant map.

\epr

\begin{proof} We may assume that the order of $\A$ is a power of  $\ell$ for some prime $\ell \neq p$.
 The arithmetic genus of $\ov Y$ equals $g$. From our assumptions it follows that the normalisation of $\ov Y$ is the projective line and that the inverse image of $P$ in the normalisation consists of one point.
Therefore $\ov Y-P$ is isomorphic to the affine line. It is well-known that $\mathbb A_{\ov F}^1$ has no $\ell$ coverings, cf. \cite[Cor. XIII. 2.12]{SGA 1}. The result now follows from Proposition \ref{simpl.con.l} applied to $\X-P$.

\end{proof}

 Let us explain the potential usefuleness of the next result, which is basically a reformulation of Proposition \ref{simpl.con.l}, convenient for applications. Calculating fundamental groups of open subvarieties is usually a difficult task. However a recent result in \cite{Schw1} allows us to do this in the case of the smooth locus for a wide class of varieties, see the remark below. In particular some of the hypothesis in the following Proposition will be automatically satisfied.

\bpr

Suppose that $\X/\O_K$ is proper but not necessarily regular, $X/K$ is  smooth, geometrically integral and $Y$ is geometrically integral.
 Assume that the maximal pro-$\ell$ quotient of $\pi_1(Y_{\ov F}-\rm{Sing}(Y_{\ov F}))$ is trivial.

 Let $\mathcal A \in \Br(X)[\ell^{\infty}]$ and assume one of the following:
 \begin{enumerate}

\item $\ell \neq p$;
\item  $\ell=p$  and $\sw_{\mathfrak p}(\A)=0$ where $\mathfrak p\in \X$ is the generic point of $Y$.

\end{enumerate}

Let $\X'=\X-\rm{Sing}(Y)$.
Then

$$
ev_{\mathcal A}: \X'(\O_K)\to \Br(K)
$$
is a constant map.

\epr

\begin{proof}
This follows from  Proposition \ref{simpl.con.l} applied to $\X'/\O_K$.

\end{proof}

\begin{remark} 1) Globally $F$-regular varieties were introduced in \cite{Smith}, and they include many classes of varieties. The following result is related to the assumptions of the previous Proposition. Suppose that $\ov Y$ is projective globally $F$-regular and the codimension of the singular locus of $Y$ is at least $2$. Then by \cite[Cor. 5.7 and Cor. 4.18]{Schw1} we have that $N=|\pi_1(Y_{\ov F}-\rm{Sing}(Y_{\ov F}))|\leq \frac{1}{s(R)} $ where $s(R)$ is the $F$-signature of a section ring of $\ov Y$ with respect to an ample sheaf and  moreover $N$ is coprime to $p$ .

2) In favourable cases, we might be able to show that $\X'(\O_K)=X(K)$ by looking at the equations of $\X$.

\end{remark}

{\bf Acknowledgements.}

The author is grateful to Martin Bright and Rachel Newton for sharing their preprint \cite{BN}, and patiently explaining its contents.
The author would like to thank Alexei Skorobogatov and Anthony V\'arilly-Alvarado for useful discussions, and the anonymous referee for numerous helpful comments.


\begin{thebibliography}{}





\bibitem{AB} Auslander, M., Brumer, A.: Brauer groups of discrete valuation rings.
Nederl. Akad. Wetensch. Proc. Ser. A 71 Indag. Math. 30 1968 286-296


\bibitem{CCSurf} Barth, W.P., Hulek, K., Peters, C.A.M., Van de Ven, A.: Compact complex surfaces.
Second edition. Ergebnisse der Mathematik und ihrer Grenzgebiete. 3. Folge. A Series of Modern Surveys in Mathematics, 4. Springer-Verlag, Berlin, 2004. xii+436 pp.


\bibitem{Bomb} Bombieri, E., Mumford, D.: Enriques' classification of surfaces in char. p. III.
Invent. Math. 35 (1976), 197-232.

\bibitem{BLR} Bosch, S., L\"utkebohmert, W., Raynaud, M.: N\'eron models.
Ergebnisse der Mathematik und ihrer Grenzgebiete (3), 21. Springer-Verlag, Berlin, 1990. x+325 pp.


\bibitem{Bri.eff} Bright, M.: Efficient evaluation of the Brauer-Manin obstruction.
Math. Proc. Cambridge Philos. Soc. 142 (2007), no. 1, 13-23.

\bibitem{Br2}  Bright, M.: Bad reduction of the Brauer-Manin obstruction. J. Lond. Math. Soc. (2),
91(3):643-666, 2015.

\bibitem{BN}  Bright, M., Newton, R.: Evaluating the wild Brauer group. preprint arXiv:2009.03282

\bibitem{Schw1}  Carvajal-Rojas, J., Schwede, K., Tucker, K.:
Fundamental groups of F-regular singularities via F-signature.
Ann. Sci. \'Ec. Norm. Sup\'er. (4) 51 (2018), no. 4, 993-1016.



\bibitem{CT} Colliot-Th\'{e}l\`{e}ne, J.-L.: The Brauer-Manin obstruction for complete intersections of dimension $\geq$ 3. Appendix to B. Poonen and J.F.
Voloch, Random Diophantine equations. Progr. Math. 226 Arithmetic of higher-dimensional
algebraic varieties (Palo Alto, CA, 2002), Birkhauser, Boston, MA, 2004, pp. 175-184.

\bibitem{CT03} Colliot-Th\'{e}l\`{e}ne, J.-L.: Points rationnels sur les fibrations. Higher dimensional varieties and rational points (Budapest, 2001), 171-221,
Bolyai Soc. Math. Stud., 12, Springer, Berlin, 2003.


\bibitem{CT.San.con} Colliot-Th\'{e}l\`{e}ne, J.-L., Sansuc, J.-J.: La descente sur les vari\'et\'es rationnelles. Journ\'ees de G\'eometrie Alg\'ebrique d'Angers, Juillet 1979, Sijthoff \& Noordhoff,
    Alphen aan den Rijn, 1980, pp. 223-237.


\bibitem{CTSk} Colliot-Th\'{e}l\`{e}ne, J.-L., Skorobogatov, A.N.: Good reduction of the Brauer-Manin obstruction.
Trans. Amer. Math. Soc. 365 (2013), no. 2, 579-590.

\bibitem{CTSkbook} Colliot-Th\'{e}l\`{e}ne, J.-L., Skorobogatov, A.N.: The Brauer-Grothendieck group.
Ergebnisse der Mathematik und ihrer Grenzgebiete. 3. Folge. A Series of Modern Surveys in Mathematics, 71. Springer-Verlag, Berlin, 2021.


\bibitem{Dolg} Cossec, F., Dolgachev, I.: Enriques surfaces. I.
Progress in Mathematics, 76. Birkh\"auser Boston, Inc., Boston, MA, 1989.



\bibitem{Deb} Debarre, O.:
Vari\'et\'es rationnellement connexes (d'apr\`es T. Graber, J. Harris, J. Starr et A. J. de Jong).
S\'eminaire Bourbaki. Vol. 2001/2002.
Ast\'erisque No. 290 (2003), Exp. No. 905, ix, 243-266.

\bibitem{Del} Deligne, P., Illusie, L.: Rel\`evements modulo $p^2$ et d\'ecomposition du complexe de de Rham.
Invent. Math. 89 (1987), no. 2, 247-270.

\bibitem{EGAIV2} Grothendieck, A.: \'El\'ements de g\'eom\'etrie alg\'ebrique. IV. \'Etude locale des sch\'emas
et des morphismes de sch\'emas, Second partie, Inst. Hautes \'Etudes Sci. Publ. Math.
24 (1965).

\bibitem{Har} Hartshorne, R.: Algebraic geometry. Graduate Texts in Mathematics, No. 52. Springer-Verlag, New York-Heidelberg, 1977.


\bibitem{Ill1} Illusie, L.: Crystalline cohomology. Motives (Seattle, WA, 1991), 43-70,
Proc. Sympos. Pure Math., 55, Part 1, Amer. Math. Soc., Providence, RI, 1994.

\bibitem{Kato1} Kato, K.: Swan conductors for characters of degree one in the imperfect residue field case. Algebraic K-theory and algebraic number theory (Honolulu, HI, 1987), 101-131,
Contemp. Math., 83, Amer. Math. Soc., Providence, RI, 1989.

\bibitem{Kato2} Kato, K.: Galois cohomology of complete discrete valuation fields. Algebraic K-theory, Part II (Oberwolfach, 1980), pp. 215-238,
Lecture Notes in Math., 967, Springer, Berlin-New York, 1982

\bibitem{Kollar} Koll\'{a}r, J.: Rational curves on algebraic varieties.
Ergebnisse der Mathematik und ihrer Grenzgebiete. 3. Folge. , 32. Springer-Verlag, Berlin, 1996.


\bibitem{Laz} Lazarsfeld, R.: Positivity in algebraic geometry. I.
Classical setting: line bundles and linear series. Ergebnisse der Mathematik und ihrer Grenzgebiete. 3. Folge. A Series of Modern Surveys in Mathematics, 48. Springer-Verlag, Berlin, 2004.



\bibitem{ShRu} Rudakov, A. N., \u Safarevi\u c, I.R.:
Inseparable morphisms of algebraic surfaces. Izv. Akad. Nauk SSSR Ser. Mat. 40 (1976), no. 6, 1269-1307, 1439.

\bibitem{Lang.Nyg} Lang, W.E., Nygaard, N.O.: A short proof of the Rudakov-\u Safarevi\u c theorem.
Math. Ann. 251 (1980), no. 2, 171-173.

\bibitem{Man} Manin, Y. I.: Le groupe de Brauer-Grothendieck en g\'eom\'etrie diophantienne. Actes du Congr\`es International des Math\'ematiciens (Nice, 1970), Tome 1, pp. 401-411. Gauthier-Villars, Paris, 1971.


\bibitem{Mats} Matsumoto, Y.: On good reduction of some K3 surfaces related to abelian surfaces.
Tohoku Math. J. (2) 67 (2015), no. 1, 83-104.


\bibitem{Milne} Milne, J.S.: \'Etale cohomology.
Princeton Mathematical Series, 33. Princeton University Press, Princeton, N.J., 1980. xiii+323 pp.


\bibitem{Poonen} Poonen, B.: Insufficiency of the Brauer-Manin obstruction applied to etale covers.
Ann. of Math. (2) 171 (2010), no. 3, 2157-2169.


\bibitem{SGA 1} Grothendieck, A.: Rev\^ etements \'etales et groupe fondamental (SGA 1).
S\'eminaire de g\'eom\'etrie alg\'ebrique du Bois Marie 1960-61. Directed by A. Grothendieck. With two papers by M. Raynaud. Updated and annotated reprint of the 1971 original [Lecture Notes in Math., 224, Springer, Berlin]. Documents Mathematiques (Paris), 3. Soci\'et\'e Math\'ematique de France, Paris, 2003.


\bibitem{Skcounter} Skorobogatov, A.N.: Beyond the Manin obstruction.
Invent. Math. 135 (1999), no. 2, 399-424

\bibitem{Skbook} Skorobogatov, A.N.:
Torsors and rational points.
Cambridge Tracts in Mathematics, 144. Cambridge University Press, Cambridge, 2001.

\bibitem{Skconj} Skorobogatov, A.N.: Diagonal quartic surfaces. Oberwolfach Rep. 33 (2009): 76-9.

\bibitem{Smith} Smith, K.E.: Globally F-regular varieties: applications to vanishing theorems for quotients of Fano varieties.
Dedicated to William Fulton on the occasion of his 60th birthday.
Michigan Math. J. 48 (2000), 553-572.

\bibitem{Span} Spanier, E.: The homology of Kummer manifolds.
Proc. Amer. Math. Soc. 7 (1956), 155-160.


\bibitem{Stacks} The Stacks Project Authors: Stacks Project \url{https://stacks.math.columbia.edu}, 2021.

\bibitem{Valla} Valla, G.: Certain graded algebras are always Cohen-Macaulay.
J. Algebra 42 (1976), no. 2, 537-548.

\bibitem{Wittrc} Wittenberg, O.: Rational points and zero-cycles on rationally connected varieties over number fields. Algebraic geometry: Salt Lake City 2015, 597-635,
Proc. Sympos. Pure Math., 97.2, Amer. Math. Soc., Providence, RI, 2018.


\end{thebibliography}


\noindent Department of Mathematics and Statistics, University of Cyprus, P.O. Box 20537,
1678, Nicosia, Cyprus

\noindent ieronymou.evis@ucy.ac.cy

\end{document}